\documentclass[12pt]{amsart}
\usepackage[utf8]{inputenc}
\usepackage[total={6.5in,9in},top=1in, left=1in, right=1in, bottom=1in]{geometry}
\usepackage{epsfig}
\usepackage{fancyhdr}
\usepackage{libertine}
\usepackage[libertine]{newtxmath}
\usepackage{color,amsmath,amssymb,mathtools,graphicx}
\usepackage{bbold}
\usepackage{fancyhdr,xcolor}
\usepackage{graphicx,amsmath,mathtools,amssymb,float}
\usepackage{bm}
\usepackage{epsfig,color,fancyhdr,setspace}
\usepackage[abs]{overpic}
\usepackage{import}
\usepackage{lineno}
\usepackage[allcolors = blue,colorlinks]{hyperref}
\usepackage[center]{caption}
\usepackage{subcaption}
\usepackage{breqn}

\newtheorem{theorem}{Theorem}[section]
\newtheorem{lemma}[theorem]{Lemma}
\newtheorem{definition}[theorem]{Definition}

\newtheorem{remark}[theorem]{Remark}

\newtheorem{conjecture}[theorem]{Conjecture}

\title{Kauffman Bracket Skein Module of the Connected Sum of Handlebodies: A Counterexample 
}

\author{Rhea Palak Bakshi}
\address{Department of Mathematics, The George Washington University, Washington DC, USA}
\email{rhea\_palak@gwu.edu}

\author{J\'{o}zef H. Przytycki}
\address{Department of Mathematics, The George Washington University, Washington DC, USA and \linebreak Department of Mathematics, University of Gda\'{n}sk, Poland}
\email{przytyck@gwu.edu}

\date{\today}

\begin{document}

\begin{abstract}

In this paper we disprove a twenty-two year old theorem about the structure of the Kauffman bracket skein module of the connected sum of two handlebodies. We achieve this by analysing handle slidings on compressing discs in a handlebody. We  find more relations than previously predicted for the Kauffman bracket skein module of the connected sum of handlebodies, when one of them is not a solid torus. Additionally, we speculate on the structure of the Kauffman bracket skein module of the connected sum of two solid tori.

\end{abstract}

\maketitle


\section{Introduction}

Skein modules were introduced by the second author in \cite{smof3} with the goal of building an algebraic topology based on knots. They generalise the skein theory of the various link polynomials in $S^3$, for example, the Alexander, Jones, Kauffman bracket, and HOMFLYPT polynomial link invariants, to arbitrary $3$-manifolds. The skein module based on the Kauffman bracket skein relation is the most comprehensively studied and best understood skein module of all. \\

Let $M$ be an oriented $3$-manifold, $\mathcal{L}^{fr}$ the set of unoriented framed links (including the empty link $\varnothing$) in $M$ up to ambient isotopy, $R$ a commutative ring with unity, and $A$ a fixed invertible element in $R$.
Consider the submodule $S_{2, \infty}^{sub}$ of the free $R$-module $R\mathcal{L}^{fr}$ generated by the Kauffman bracket skein relation, $L_+ - AL_0 - A^{-1}L_{\infty}$, and
 the trivial component relation, $L \sqcup {\pmb \bigcirc}  + (A^2 + A^{-2})L$, where $\pmb{\bigcirc}$ denotes the trivial framed knot and the skein triple $(L_+$, $L_0$, $L_{\infty})$ denotes three framed links in $M$ which are identical except in a small $3$-ball in $M$ where they differ as shown: 

\[  \begin{minipage}{1.2 in} \includegraphics[width=\textwidth]{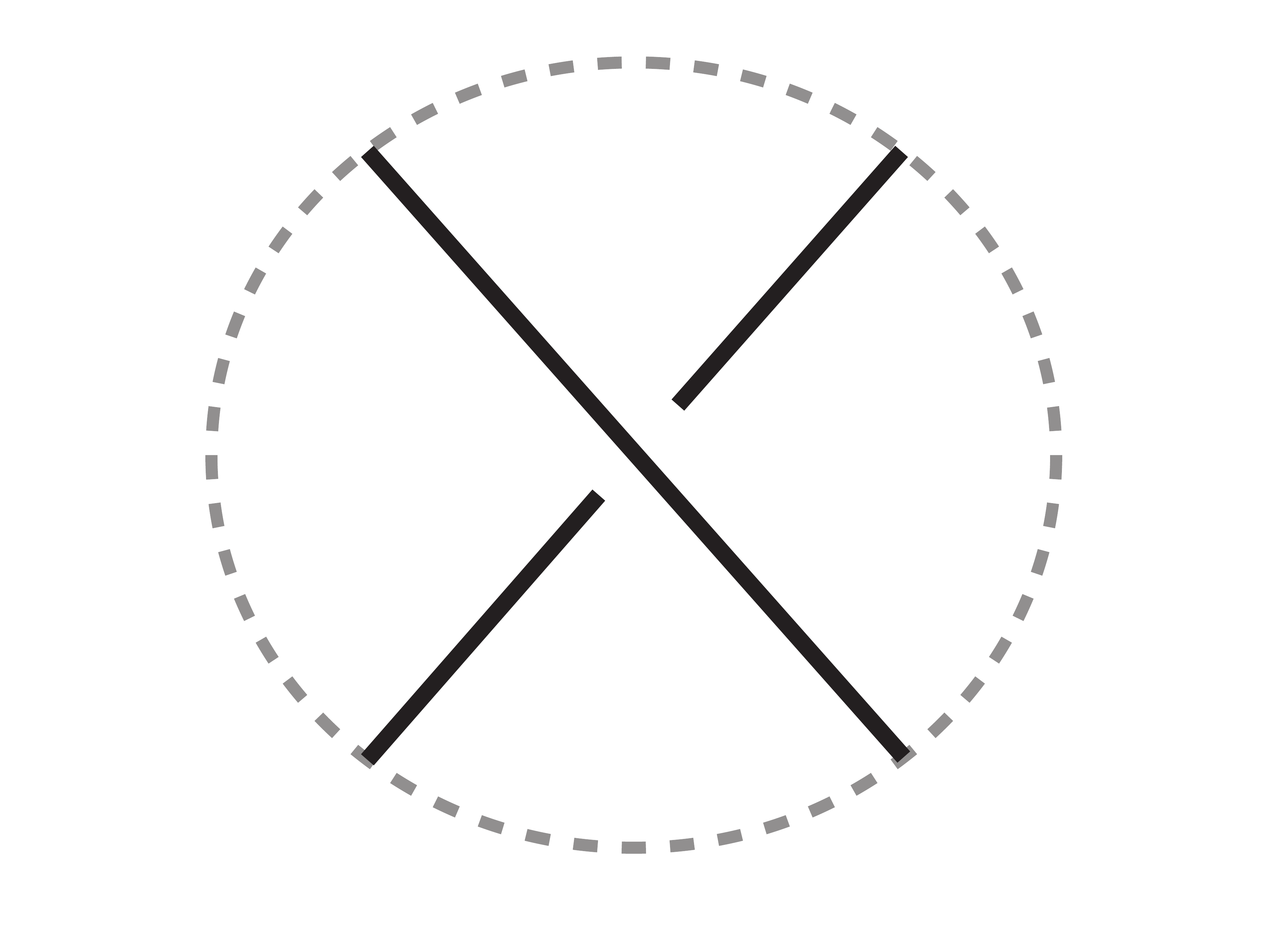} \vspace{-15pt} \[L_+\] \end{minipage} 
               \qquad
        \begin{minipage}{1.2 in}\includegraphics[width=\textwidth]{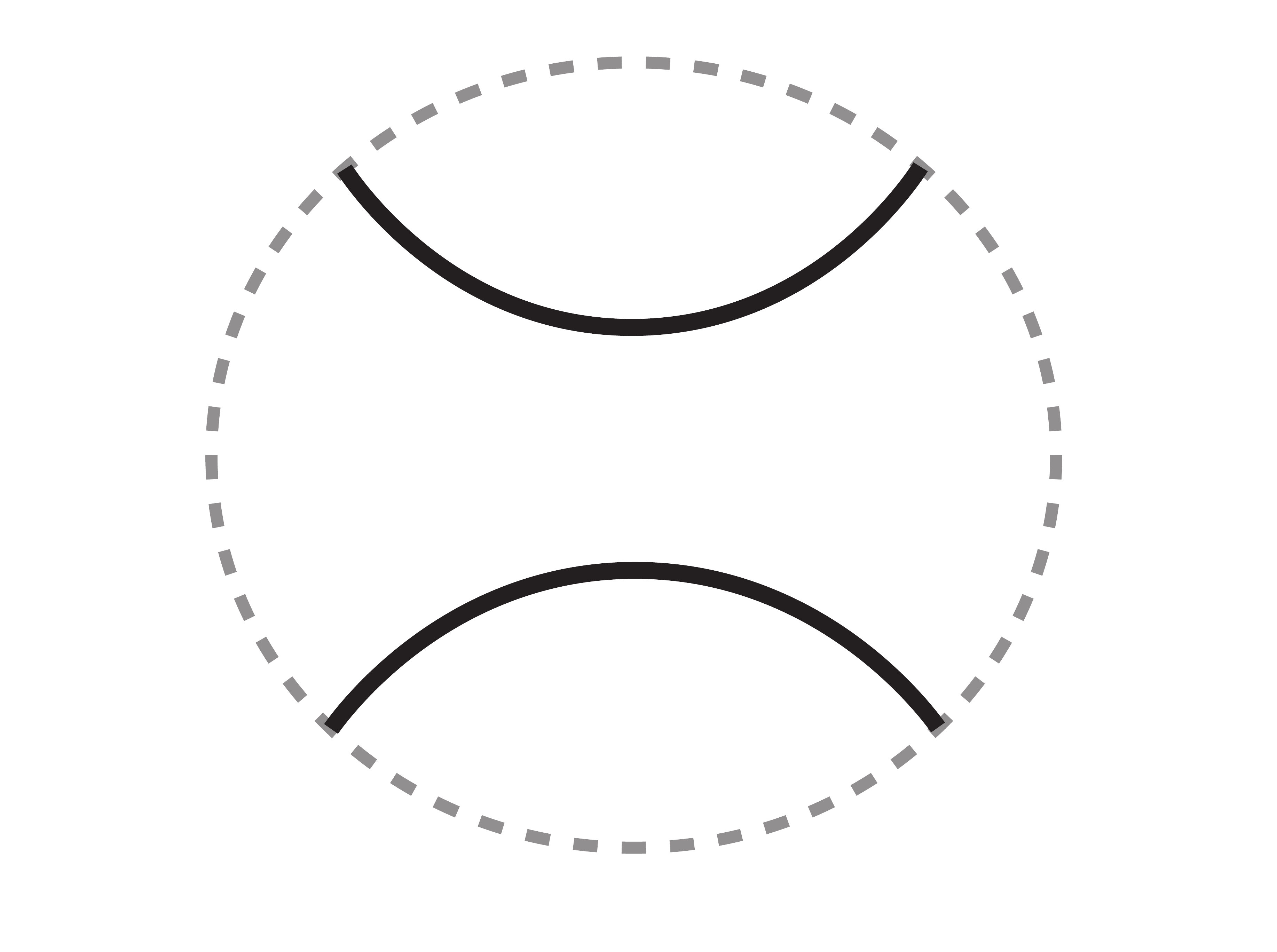} \vspace{-15pt} \[L_0\] \end{minipage}
         \qquad
        \begin{minipage}{1.2 in}\includegraphics[width=\textwidth]{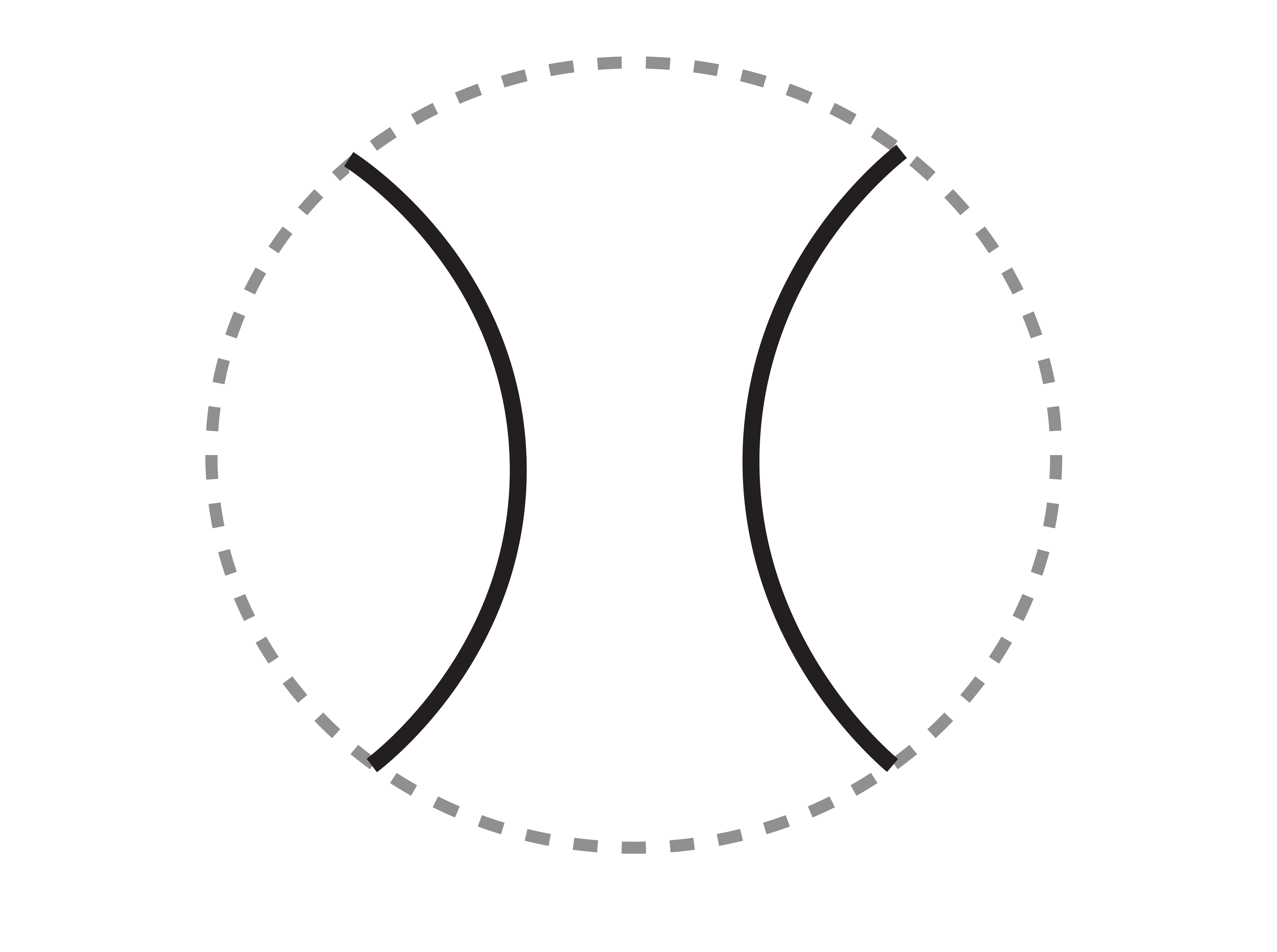} \vspace{-15pt} \[ L_\infty\]\end{minipage} 
        \]

The {\bf Kauffman bracket skein module} (KBSM) of $M$ is defined as the quotient $\mathcal{S}_{2,\infty}(M;R,A) = R\mathcal{L}^{fr}/S_{2, \infty}^{sub}.$ For brevity, when $R = \mathbb{Z}[A^{\pm 1}]$ we use the notation $\mathcal{S}_{2,\infty}(M)$  in the remainder of the paper. One can also work with the relative case in which the oriented $3$-manifold $(M, \partial M)$ has $2k$ marked points $\{x_1, x_2, \hdots, x_{2k}\}$ on $\partial M$. 
The {\bf relative Kauffman bracket skein module} (RKBSM) of $(M, \partial M, \{x_i\}_1^{2k})$, denoted by $\mathcal{S}_{2,\infty}(M, \{x_i\}_1^{2k}; R, A)$, is the set of all ambient isotopy classes of relative framed links $L$ in $(M, \partial M, \{x_i\}_1^{2k})$ keeping $\partial M$ fixed, such that $L \cap \partial M = \partial L = \{x_i\}_1^{2k}$,  modulo the Kauffman bracket relations (see \cite{fundamentals}).

In the year 2000, the second author published the following fundamental theorem about the connected sum of the Kauffman bracket skein module of oriented $3$-manifolds over the ring $\mathbb{Z}[A^{\pm 1}]$ localised by inverting all the cyclotomic polynomials in $A$.
 
 \begin{theorem}\cite{connsum}\label{connsum}

If $M$ and $N$ are compact, oriented $3$-manifolds, $M \ \# \ N$ denotes their connected sum, and $A^k - 1$ is invertible in $R$ for any $k >0$, then $$\mathcal{S}_{2,\infty}(M \ \# \  N; R, A) = \mathcal{S}_{2,\infty}(M; R, A) \otimes \mathcal{S}_{2,\infty}(N; R, A).$$

\end{theorem}

 In particular, the result holds when $R = \mathbb{Q}(A)$. This theorem was an an essential tool used in \cite{wittenresolved} to resolve Witten's finiteness conjecture for the KBSM of $3$-manifolds over $\mathbb{Q}(A)$, since it shows that if $\mathcal{S}_{2,\infty}(M; \mathbb{Q}(A))$ and $\mathcal{S}_{2,\infty}(N; \mathbb{Q}(A))$ are finite dimensional, then so is $\mathcal{S}_{2,\infty}(M \ \# \ N; \mathbb{Q}(A))$, that is, finiteness of KBSMs is stable under connected sums. 
 However, Theorem \ref{connsum} does not hold when the ring $R = \mathbb{Z}[A^{\pm 1}]$. In this case $\mathcal{S}_{2,\infty}(M)$ is not always finitely generated and often contains torsion (see \cite{s1s2}). In fact, much less is known about the structure of the Kauffman bracket skein module of an oriented $3$-manifold when $R = \mathbb{Z}[A^{\pm 1}]$ than when $R = \mathbb{Q}(A)$. Julien March\'e had proposed a conjecture (see \cite{basissurfacetimess1}) about the structure of the KBSM over $R = \mathbb{Z}[A^{\pm 1}]$ which was recently disproved by the first author in \cite{rheasolomarche}. \\
 
 With the goal of understanding the structure of Kauffman bracket skein modules over the ring $\mathbb{Z}[A^{\pm 1}]$ and giving a complete and detailed description of the KBSM of connected sums and disc sums, the second author had stated the following theorem without proof (Theorem 7.1 in \cite{connsum}) about the Kauffman bracket skein module over $\mathbb{Z}[A^{\pm 1}]$ of the connected sum of handlebodies.

\begin{figure}[h]
    \centering
    \includegraphics[scale=0.6]{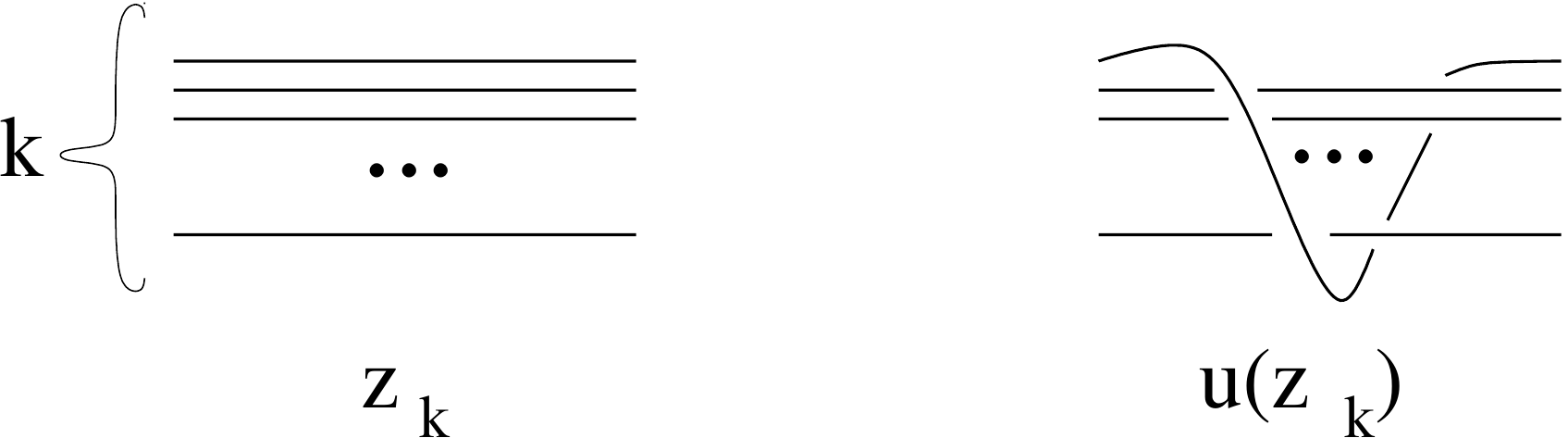}
    \caption{{\color{white}.}}
    \label{fig7.1}
\end{figure}

\begin{figure}[h]
    \centering
    \includegraphics[scale=0.6]{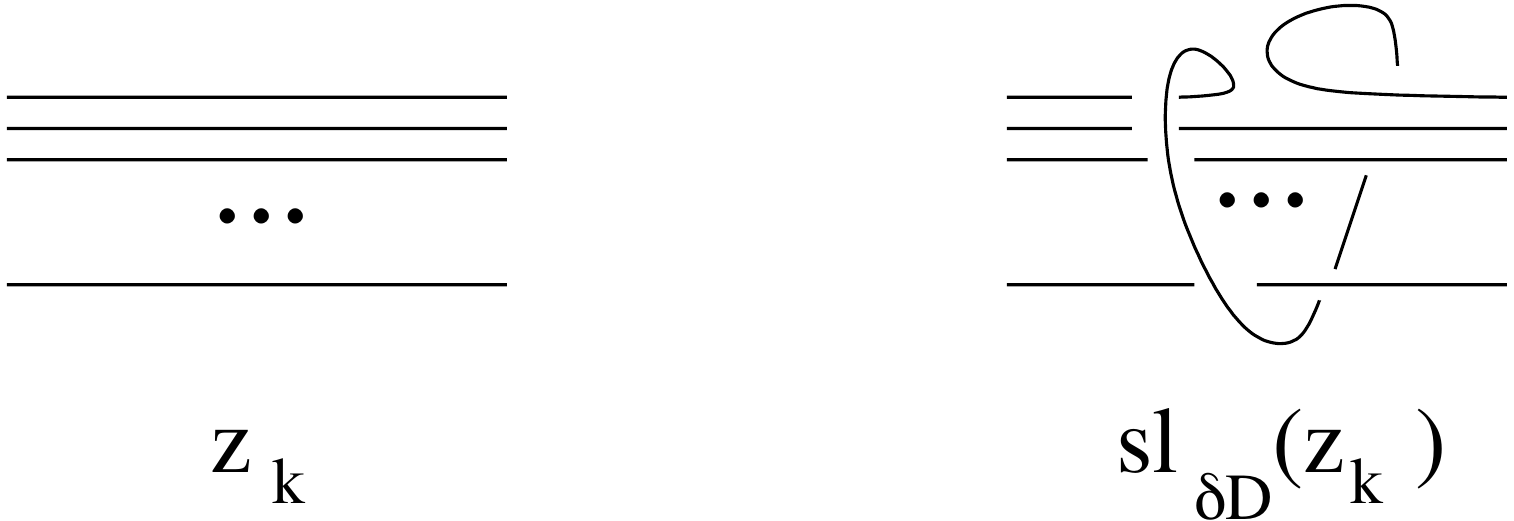}
    \caption{{\color{white}.}}
    \label{fig7.2}
\end{figure}

\begin{theorem}\cite{connsum}\label{przytyckitheorem}

Let $H_n$ denote a genus $n$ handlebody and $F_{0,n+1}$ be a disc with $n$ holes so that $H_n = F_{0,n+1} \times I$. Then,

$${\mathcal S_{2,\infty}}(H_n \ \# \ H_m) = {\mathcal S_{2,\infty}}(H_{n+m})/\mathcal I,$$
where $\mathcal I$ is the ideal generated by expressions $z_k-A^6u(z_k)$, 
for any even $k\geq 2$, and $z_k \in B_k(F_{0,n+m+1})$, where 
$B_k(F_{0,n+m+1})$ is composed of links without contractible components and
with geometric intersection number $k$ with a disc $D$ separating
$H_n$ and $H_m$. $u(z_k)$ is a modification of $z_k$
in the neighbourhood of $D$, as shown in Figure \ref{fig7.1}.
The relation $z_k = A^6u(z_k)$, is a result of the sliding relation
$z_k = sl_{\partial D}(z_k)$ as illustrated  in Figure \ref{fig7.2}.

\end{theorem}

In this paper we disprove this theorem by providing a counterexample which is given by $H_n \ \# \ H_m$, $n \geq 2, m \geq 1$ and we show that the ideal $\mathcal I$ should be replaced by a strictly bigger ideal to obtain the equality in Theorem \ref{przytyckitheorem}. At the time of writing this paper the case of $H_1 \ \# \ H_1$ is still open (see Conjecture \ref{h1h1conjecture}).

\section{Preliminaries}

In this section we discuss several important properties of Kauffman bracket skein modules that are pertinent to this paper, including a description of the KBSM of any $3$-manifold using generators and relations.

\begin{theorem}\cite{smof3, fundamentals}\label{propkbsm}

\begin{enumerate}

    \item Let $i : M \hookrightarrow N$ be an orientation preserving embedding of $3$-manifolds. This yields a homomorphism of skein modules $i_* : \mathcal{S}_{2,\infty}(M;R,A) \longrightarrow
\mathcal{S}_{2,\infty}(N;R,A)$. This correspondence leads to a functor from the
category of $3$-manifolds and orientation preserving embeddings (up to ambient isotopy) to the category of $R$-modules with a specified invertible
element $A \in R$.\label{embhom} \\

    \item  If $N$ is obtained from $M$ by adding a $3$-handle 
     to $M$ and $i : M \hookrightarrow N$ is the associated embedding, then
$i_* : \mathcal{S}_{2,\infty}(M;R,A) \longrightarrow \mathcal{S}_{2,\infty}(N;R,A)$ is an isomorphism.\label{3hand} \\

\item \label{disjointunion} If $M_1 \sqcup M_2$ is the disjoint sum of oriented $3$-manifolds $M_1$ and $M_2$ then
$\mathcal{S}_{2,\infty}(M_1 \sqcup M_2;R,A) = \mathcal{S}_{2,\infty}(M_1;R,A) \otimes  \mathcal{S}_{2,\infty}(M_2;R,A)$.\\

\item {\it (The Universal Coefficient Property)}

Let $R$ and $R'$ be commutative rings with unity and $r : R \longrightarrow R'$ be a homomorphism. 
Then the identity map on $\mathcal{L}^{fr}$
induces the following isomorphism of $R'$ (and R) modules:
$$\overline{r} : \mathcal{S}_{2,\infty}(M;R,A) \otimes_R  R' \longrightarrow \mathcal{S}_{2,\infty}(M;R', r(A)).$$

\end{enumerate}

\end{theorem}

The following lemma allows one to write a presentation of the Kauffman bracket skein module of any compact oriented $3$-manifold using its Heegaard decomposition and knowledge of the presentation of the KBSM of any handlebody. Theorem \ref{surfacecrossint} describes the KBSM and RKBSM of trivial surface $I$-bundles and in particular, handlebodies.

\begin{lemma}[Handle Sliding Lemma]\cite{fundamentals,connsum, kbsmlens}

\begin{enumerate}

\item Let $M$ be a $3$-manifold with boundary $\partial M$ and $\gamma$ be a
simple closed curve on $\partial M$. Let $N = M_{\gamma}$ be the $3$-manifold
obtained from $M$ by adding a $2$-handle along $\gamma$ and $i : M \hookrightarrow N$ be the associated embedding. Then
$i_* : \mathcal{S}_{2,\infty}(M;R,A) \longrightarrow \mathcal{S}_{2,\infty}(N;R,A)$ is an epimorphism.\label{2hand} Furthermore, the kernel of $i_*$ is generated by the relations yielded by $2$-handle sliding. More precisely, if $\mathcal{L}^{fr}_{gen}$
is a set of framed links in $M$ which generates $\mathcal{S}_{2,\infty}(M;R,A)$,
then $\mathcal{S}_{2,\infty}(N;R,A) = \mathcal{S}_{2,\infty}(M;R,A)/J$, where $J$ is the submodule of
$\mathcal{S}_{2,\infty}(M;R,A)$ generated by the expressions $L - sl_{\gamma}(L)$.  Here $L \in \mathcal{L}^{fr}_{gen}$ and $sl_{\gamma}(L)$ is obtained from $L$ by sliding it along $\gamma$ (that is, we perform $2$-handle sliding).\\

\item Let $M$ be an oriented compact $3$-manifold and consider its Heegaard decomposition
(that is, M is obtained from the handlebody $H_m$ by adding
$2$ and $3$-handles to it). Then the generators
of $\mathcal{S}_{2,\infty}(M;R,A)$ are generators of $\mathcal{S}_{2,\infty}(H_m;R,A)$ and the relators of $\mathcal{S}_{2,\infty}(M;R,A)$
are yielded by $2$-handle slidings.

\end{enumerate}

\end{lemma}

\begin{theorem}\cite{smof3,fundamentals}\label{surfacecrossint}

\begin{enumerate}

    \item $\mathcal{S}_{2,\infty}(F \times [0,1])$ is a free module generated by the empty link $\varnothing$ and links in $F$ which have no trivial components. Here $F$ is an oriented surface and each link in $F$ is equipped with an arbitrary, but specific framing. This applies in particular to handlebodies, since $H_{n} = F_{0,n+1} \times I$, where $H_n$ is a handlebody of genus $n$ and $F_{g,b}$ denotes a surface of genus $g$ with $b$ boundary components. \\
    
    \item If $\partial F \neq \varnothing$ and all $x_i, 1 \leq i \leq 2n$, lie on $\partial F \times \{\frac{1}{2}\}$, then $\mathcal S_{2, \infty}(F \times I, \{x_i\}_1^{2n}; R, A)$ is a free $R$-module whose basis is composed of relative links in $F$ without trivial components.
 
 \end{enumerate}

\end{theorem}

The following definition will be useful in our construction of the counterexample to Theorem \ref{przytyckitheorem}.

\begin{definition}\label{bluedefinition}

Consider two oriented $3$-manifolds $(M, \partial M, \{x_i\}_1^{2k})$ and $(N, \partial N, \{y_j\}_1^{2k})$. Let $F_M$ be a surface in $\partial M$ containing $\{x_i\}_1^{2k}$,  $F_N$ be a surface in $\partial N$ containing $\{y_j\}_1^{2k}$ and consider the homeomorphism $f: F_M \longrightarrow F_N$ such that $\{x_i\}_1^{2k} \mapsto \{y_j\}_1^{2k}$. If $W$ is the $3$-manifold obtained by gluing $M$ and $N$ along a part of their boundaries using $f$ then we have the following bilinear form: $$\langle \ . \ \rangle : \mathcal{S}_{2,\infty}(M, \{x_i\}_1^{2k}; R, A) \times \mathcal{S}_{2,\infty}(N, \{y_j\}_1^{2k}; R, A) \longrightarrow \mathcal{S}_{2,\infty}(W; R, A).$$

\end{definition}

\section{Handle Sliding Relations in Handlebodies}

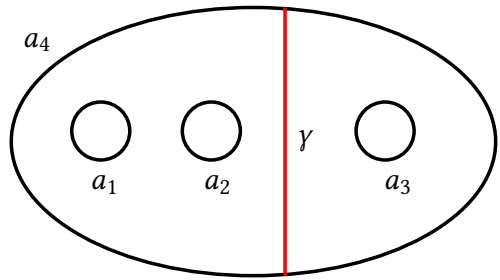
\begin{figure}[h]
     \centering
$$\vcenter{\hbox{
  \begin{overpic}[unit=1mm, scale = 0.7]{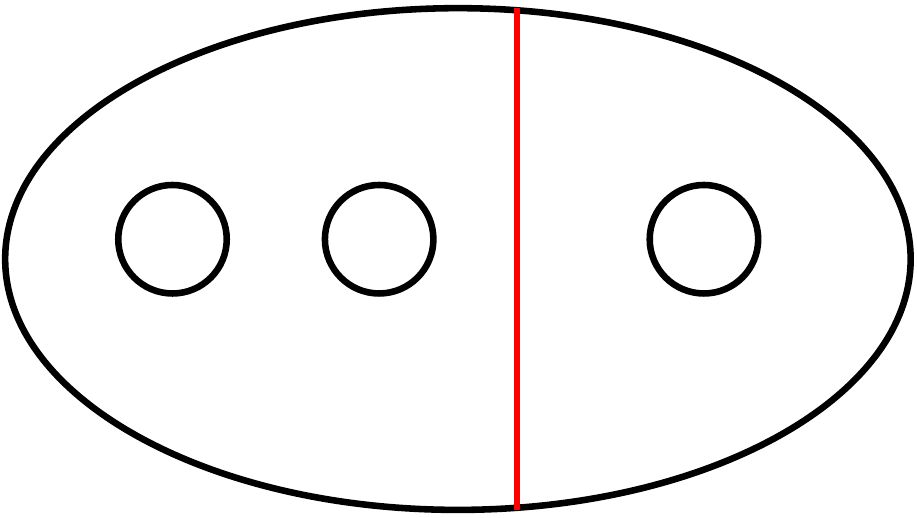} 
 \put(11,12){$a_1$}
  \put(26,12){$a_2$}
   \put(50,12){$a_3$}
    \put(2,31){$a_4$}
 \put(38.5,18){$\gamma$}
\end{overpic}}}$$
\caption{Projection of $H_2 \ \# \ H_1$ onto $F_{0,4}$}
     \label{h2connh1nocurves}
 \end{figure}

  Consider the oriented $3$-manifold $M = H_n \  \# \ H_m$, $n,m \geq 1$ and let $D$ be the compressing disc in $M$ which separates $H_n$ and $H_m$. In addition, let $\gamma = \partial D$.  Now $H_n \  \# \ H_m$ is homeomorphic to $(H_{n+m})_{\gamma}$ which is the $3$-manifold $H_{n+m}$ with a $2$-handle added along $\gamma$. Let $F_{0,n+m+1}$ be a sphere with $n+m+1$ boundary components denoted by $a_1, a_2, \hdots , a_{n+m+1}$. Then $H_{n+m} \cong F_{0,n+m+1} \times I$ and we project links in $(H_{n+m})_{\gamma}$ onto $F_{0,{n+m+1}}$ and work with their corresponding link diagrams. The projection of $\gamma$ onto $F_{0,n+m+1}$ is represented by a red line segment as illustrated in Figure \ref{h2connh1nocurves}. \\

\begin{figure}[h]
      \centering
     \begin{subfigure}{0.49 \textwidth}
     \centering$$\vcenter{\hbox{
     \begin{overpic}[unit=1mm, scale = 0.2]{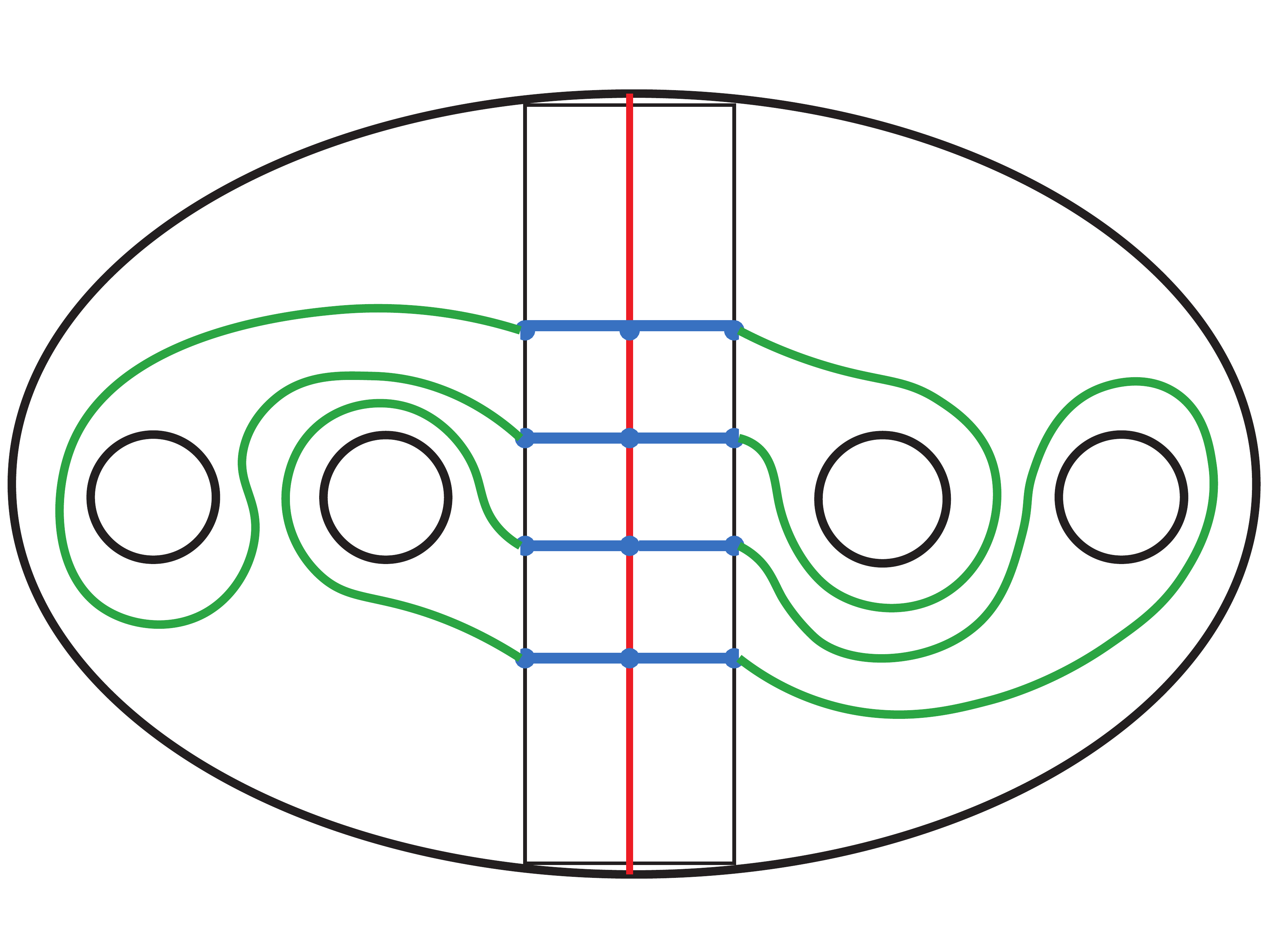} \put(7,25.5){$a_1$}
     \put(20.5,25.5){$a_2$}
     \put(50,25.5){$a_3$}
     \put(63.5,25.5){$a_4$}
     \put(8,46){$a_5$}
     \end{overpic}}}
     $$
      \caption{ $\kappa \cup Id_4$ }
      \label{z4figure}
      \end{subfigure}
      \begin{subfigure}{0.49 \textwidth}
     \centering
     $$\vcenter{\hbox{
     \begin{overpic}[unit=1mm, scale = 0.2]{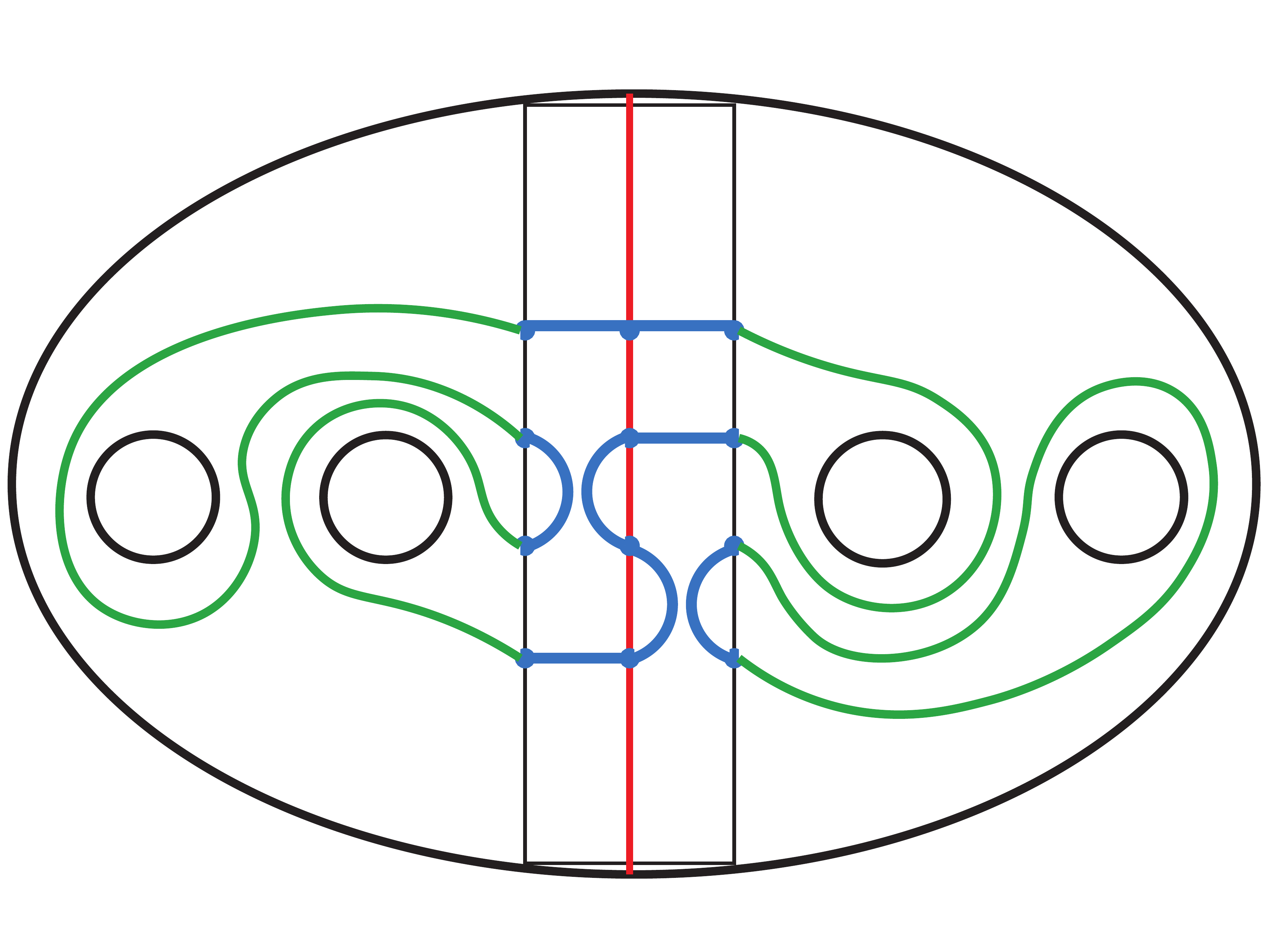} 
     \put(7,25.5){$a_1$}
     \put(20.5,25.5){$a_2$}
     \put(50,25.5){$a_3$}
     \put(63.5,25.5){$a_4$}
     \put(62,46){$a_5$}
     \end{overpic}}}$$
      \caption{$\kappa \cup e_2e_3$}
      \end{subfigure}
      \caption{Multicurves in $F_{0,5}$ with $k=4$}
      \label{connsumhandlebodytlnembdedding}
  \end{figure}
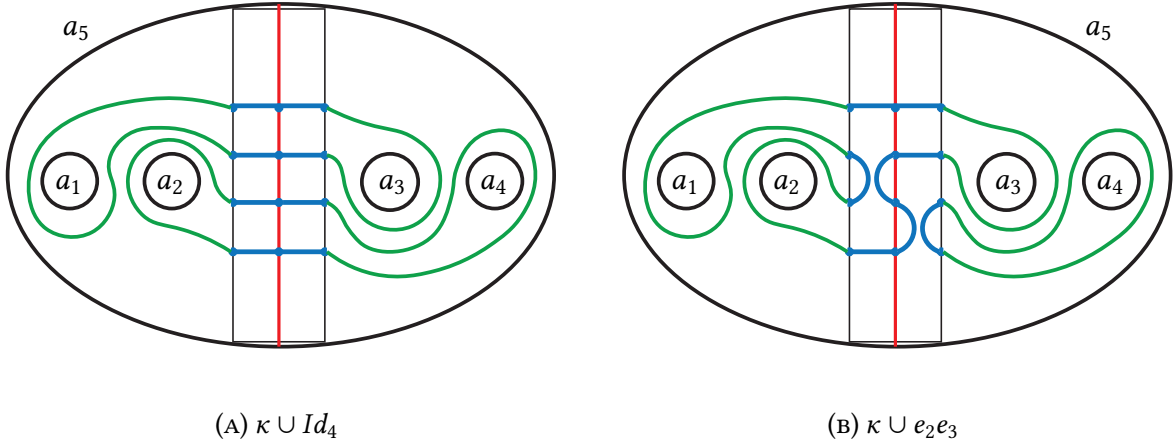

  As a special case of Definition \ref{bluedefinition}, let $W = F_{0,n+m+1} \times I$. Consider a rectangle $\mathcal{R}$, which is the regular neighbourhood of the red line segment, and its embedding under $\rho$ into $F_{0,n+m+1} = F_{0,n+m+1} \times \{\frac{1}{2}\}$. Choose $2k$ marked points on the boundary of this rectangle, with $k$ points on the left edge and $k$ points on the right edge. Consider the Temperley-Lieb box $M = \mathcal{R} \times I$ and relative links in $M$  modulo the Kauffman bracket relations. This gives rise to the Temperley-Lieb module $TL_k$ which is the relative Kauffman bracket skein module of $\mathcal{R} \times  I$.\footnote{We use Kauffman's diagrammatic visualisation \cite{tlnkauffman} of the Temperley-Lieb algebra $TL_k$, where $Id_k$ denotes the identity element and $e_i$ denotes caps connecting the marked points $x_i$ with $x_{i+1}$ and $x_{2k-i}$ with $x_{2k-i+1}$.} For any relative multicurve $\kappa$ in $F_{0,n+m+1} \setminus \mathcal{R}$ we have a module homomorphism $\rho_*: \mathcal{S}_{2,\infty}(\mathcal{R} \ \times \ I, \{x_i\}_1^{2k}) \longrightarrow \mathcal{S}_{2,\infty}(H_{n+m})$ as follows: for a given Temperley-Lieb element its image under $\rho_*$ is obtained by gluing it to $\kappa$ along the $2k$ marked points. See Figure \ref{connsumhandlebodytlnembdedding} for an example. \\

\begin{figure}[h]
        \centering
$$\vcenter{\hbox{
\begin{overpic}[unit=1mm, scale = 0.8]{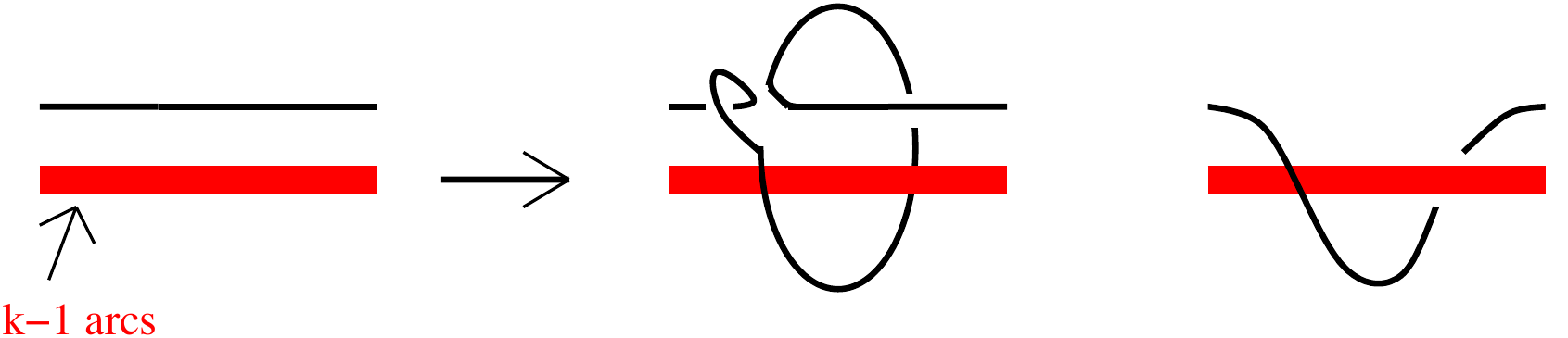}
\put(41.5,16){$\phi_t$}
\put(117,0){$u(Id_k)$}
\put(38,8){{\it sliding}}
\put(93,13){\Large{$=A^6$}}
\end{overpic}}}$$
\caption{Positive $2$-handle sliding on the upper arc}
\label{Positive handle sliding on the upper arc}
    \end{figure}

\begin{figure}[h]
    \centering
\begin{subfigure}{1.0 \textwidth}
\centering
$$\vcenter{\hbox{
\begin{overpic}[unit=1mm, scale = 0.8]{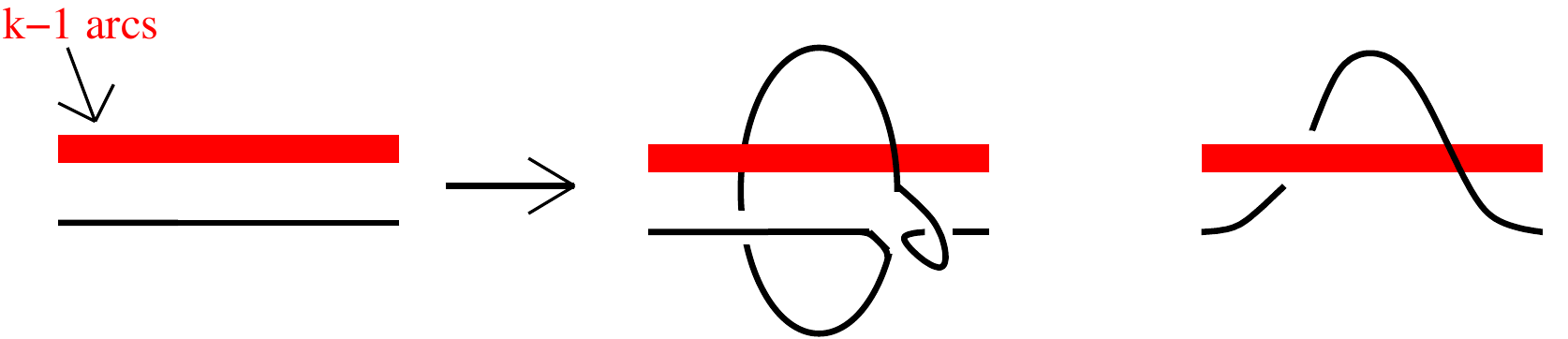}
\put(41.5,15.5){$\phi_{\ell}$}
\put(113,3){$w(Id_k)$}
\put(38,8){{\it sliding}}
\put(90,13){\Large{$=A^6$}}
\end{overpic}}}$$
\caption{Positive $2$-handle sliding on the lower arc}
\label{lowerphil}
\end{subfigure}
\begin{subfigure}{1.0 \textwidth}
\centering
$$\vcenter{\hbox{
    \begin{overpic}[unit=1mm, scale = 0.8]{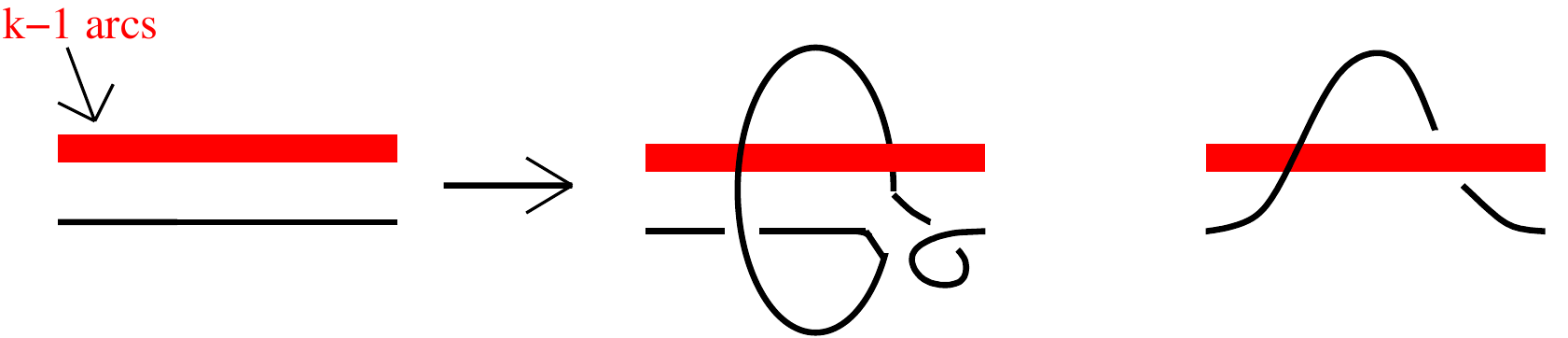}
\put(41.5,15.5){$\overline{\phi}_{\ell}$}
\put(113,3){$\overline{w}(Id_k)$}
\put(38,8){{\it sliding}}
\put(90,13){\Large{$=A^{-6}$}}
\end{overpic}}}$$
\caption{Negative $2$-handle sliding on the lower arc}
\label{lowerphilmirror}
\end{subfigure}
\caption{{\color{white}.}}
\label{Sliding on the lower arc}
\end{figure}

Let $z_k$ be a multicurve in $F_{0,n+m+1} \times \{\frac{1}{2}\}$ obtained by gluing the identity element $Id_k$ of $TL_k$ to some $\kappa$ using the homomorphism $\rho_*$ such that $z_k$ is in general position with the compressing disc $D$  having geometric crossing number $k$ with it (see Figure \ref{z4figure}).  In $(H_{n+m})_{\gamma}$ consider the $2$-handle slidings of $z_k$ along $\gamma$ described in Figures \ref{Positive handle sliding on the upper arc} and \ref{Sliding on the lower arc}. These handle slidings have  support in $(\mathcal{R} \times I)_{\gamma}$, the Temperley-Lieb box with a $2$-handle attached along $\gamma$.

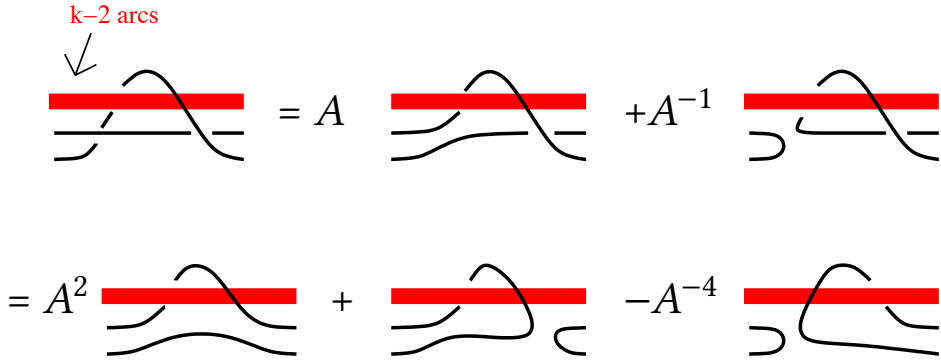
\begin{figure}[h]
    \centering
    $$\vcenter{\hbox{
    \begin{overpic}[unit=1mm, scale = 0.7]{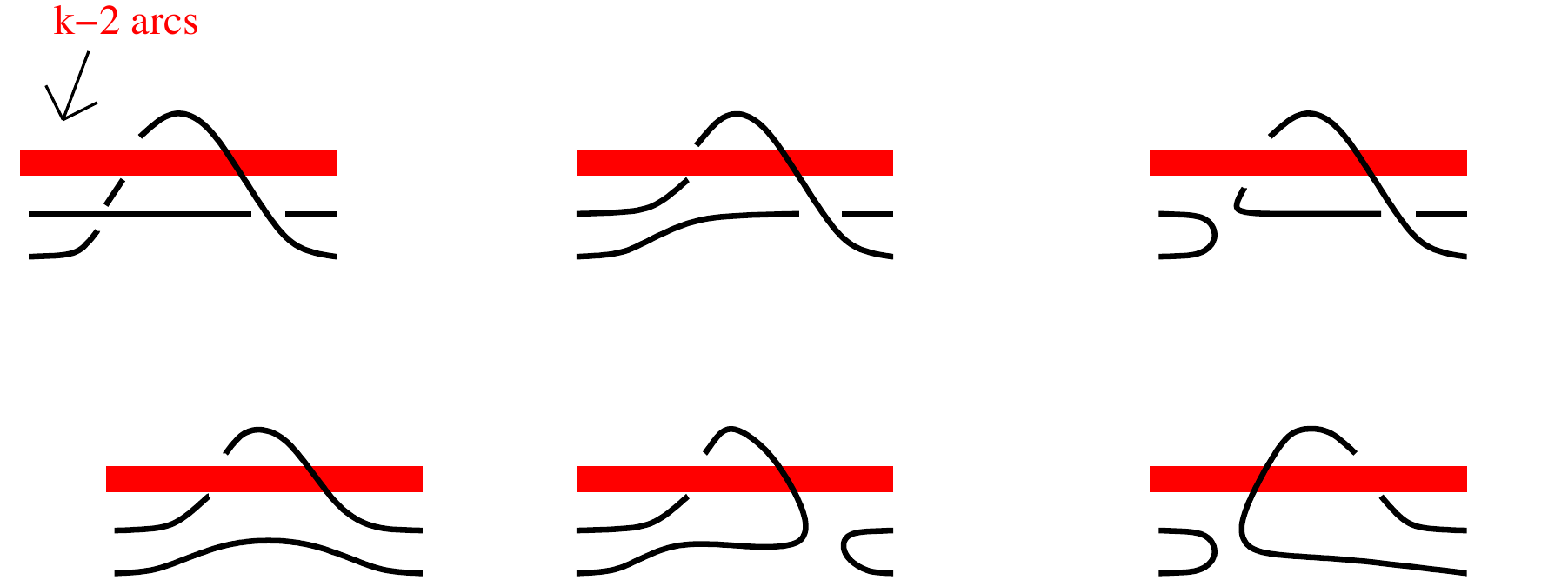}
    \put(32,31){\Large{$=A$}}
    \put(78,31){\Large{$+A^{-1}$}}
    \put(-4,6){\Large{$=A^{2}$}}
     \put(39,6){\Large{$+$}}
     \put(78,6){\Large{$-A^{-4}$}}
    \end{overpic}}}$$
    \caption{Calculation for $w(Id_k)$}
    \label{XXX}
\end{figure}

The sliding relation $Id_k=\phi_{t}(Id_k)$, coming from Figure \ref{Positive handle sliding on the upper arc}, is obtained by positive handle sliding on the top arc and the sliding relation $Id_k=\phi_{\ell}(Id_k)$, coming Figure \ref{Sliding on the lower arc}, is obtained by positive handle sliding on the lower arc. Here, $\phi_{\ell}(Id_k)= A^6w(Id_k)$, where $w(z_k)$ is the figure on the left hand side of the equation in Figure \ref{XXX}, and $\phi_t(Id_k) = A^6u(Id_k)$.  Our calculations are carried out in the Templerley-Lieb module and thus, the relations arising from the $2$-handle slidings are written in terms of the Temperley-Lieb elements. By extension, using Definition \ref{bluedefinition} and embedding $TL_{k}$ into $H_{n+m}$ as described earlier, we get the relations in $\mathcal{S}_{2,\infty}(H_{n+m})_{\gamma}$.

\section{Counterexample to Theorem \ref{przytyckitheorem}}

In this section we construct a counterexample to Theorem \ref{przytyckitheorem}. Our result is summarised as follows: 

\begin{theorem}\label{counterexampletheoremconnsum}

\begin{enumerate}

\item $\mathcal{S}_{2,\infty}(H_2 \ \# \ H_1) \neq \mathcal{S}_{2,\infty}(H_{3})/\mathcal{I}$.

\item In general, $\mathcal{S}_{2,\infty}(H_n \ \# \ H_m) \neq \mathcal{S}_{2,\infty}(H_{n+m})/\mathcal{I}$, where $n + m\geq 3$.

\end{enumerate}

\end{theorem}

 \begin{figure}[h]
    \centering
    \includegraphics[scale=0.7]{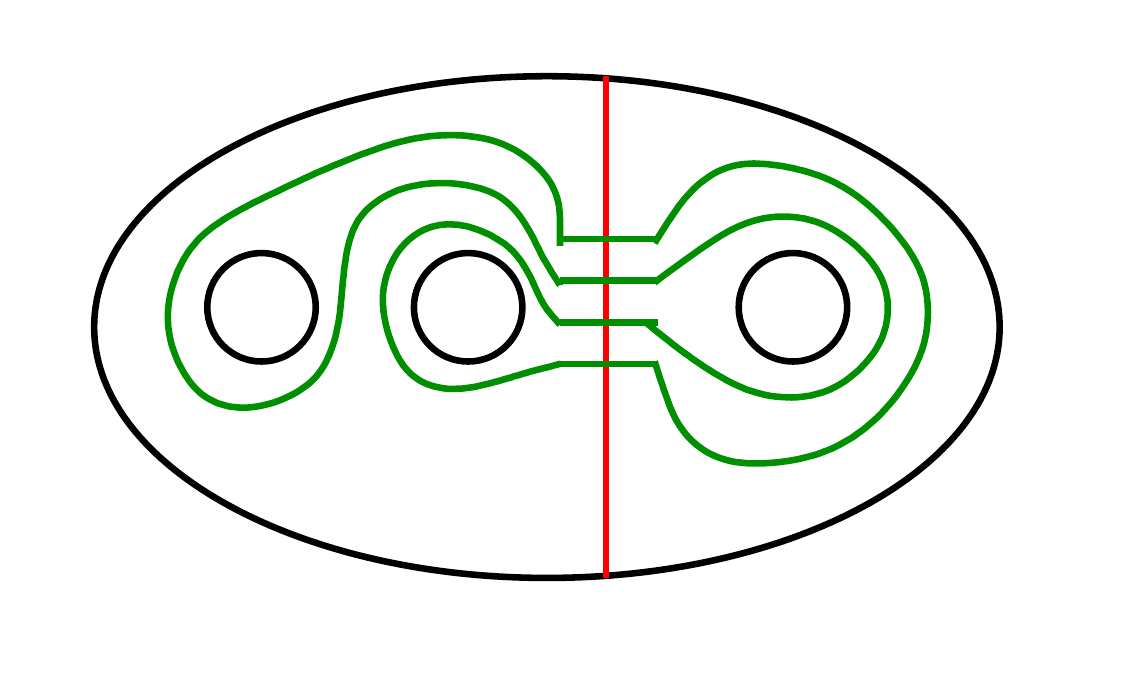}
    \caption{Curve system in $F_{0,4}$ leading to the counterexample}
    \label{H2H1counterexamplefigure}
\end{figure}

\begin{proof}

We first prove part (1) of the theorem. Part (2) follows by an easy generalisation. Consider the oriented $3$-manifold $H_2 \ \# \ H_1$ and the positive $2$-handle sliding $\phi_{\ell}$ on the bottom arc illustrated in Figure \ref{Sliding on the lower arc}. We compute $w(Id_k)$ recursively starting from $k=2$, in which case we obtain the following result.  

\begin{equation}\label{wid2}
w(Id_2)= A^2Id_2 + (1-A^{-4})e_1=A^2Id_2 +A^{-4}(A^4-1)e_1.\end{equation}

In general, from Figure  \ref{XXX}, we get the following recursive relation for $w(Id_k)$: 

\begin{equation}\label{widkrecursiverelation}
w(Id_k)=A^2w(Id_{k-1})\otimes Id_1+ (w(Id_{k-1})\otimes Id_1)\circ e_{k-1} -A^{-4}e_{k-1}\circ( \overline{w}(Id_{k-1})\otimes Id_1).
\end{equation}

Here $\overline{w}(Id_k)$ denotes the mirror image of $w(Id_k)$ (see Figure \ref{lowerphilmirror}). In particular, when $k=3$ we get the following equation: 

\begin{equation}\label{wid3equation}
w(Id_3)=A^4Id_3 +(A^2-A^{-6})e_2 + (A^2-A^{-2})e_1 + (1-A^{-4})(e_1e_2+e_2e_1),
\end{equation}

and when $k=4$, we get the following equation using Equations (\ref{widkrecursiverelation}) and (\ref{wid3equation}) :


\begin{dmath}w(Id_4)=A^6Id_4+ (A^4-1)e_1+ (A^4-A^{-4})e_2+ (A^2-A^{-2})(e_1e_2+e_2e_1)+(A^4-A^{-8})e_3 + (A^2-A^{-6})(e_1e_3 +e_2e_3+e_3e_2)+(1-A^{-4})(e_1e_2e_3+e_3e_2e_1 + e_1e_3e_2 +e_2e_3e_1).\end{dmath}

Since $\phi_{\ell}(Id_4) = A^6w(Id_4)$, therefore,

\begin{dmath}
\phi_l(Id_4) = A^{12}Id_4 + (A^{10} - A^{6})e_1 + (A^{10}-A^{-2})e_3 +  + (A^{10}- A^{2})e_2 + (A^8 - A^{4})(e_2e_1+e_1e_2)$$ $$ + (A^8-1)[e_2e_3 + e_3e_2 + e_1e_3] + (A^{6}-A^{2}) (e_3e_2e_1+e_1e_3e_2 + e_2e_3e_1+e_1e_2e_3).
\end{dmath}

Consider the handle sliding relation $\phi_l(Id_4) \equiv Id_4$ in the relative Kauffman bracket skein module of $(\mathcal{R} \times I)_\gamma$. Therefore,

\begin{dmath}\label{connsumeqn1}(1-A^{12})Id_4 \equiv (A^{10} - A^{6})e_1  + (A^{10}-A^{-2})e_3 + + (A^{10}- A^{2})e_2 + (A^8 - A^{4})(e_2e_1+e_1e_2)+ (A^8-1)[e_2e_3 + e_3e_2 + e_1e_3] + (A^{6}-A^{2}) (e_3e_2e_1+e_1e_3e_2 + e_2e_3e_1+e_1e_2e_3).\end{dmath}

By performing $2$-handle sliding $\phi_t$ on the upper string, the roles of $e_1$ and $e_3$ are exchanged in the above relation and we get the following:

\begin{dmath}\label{connsumeqn2}
(1-A^{12})Id_4 \equiv (A^{10}-A^{-2})e_1 + (A^{10} - A^{6})e_3 + (A^{10}- A^{2})e_2 + (A^8 - A^{4})(e_2e_3+e_3e_2)+ (A^8-1)(e_2e_1 + e_1e_2 + e_1e_3) + (A^{6}-A^{2}) (e_1e_2e_3 +e_1e_3e_2  + e_2e_3e_1 + e_3e_2e_1). 
\end{dmath}

Subtracting the Equation (\ref{connsumeqn1}) from Equation (\ref{connsumeqn2}) we get:

$$0 \equiv (A^6-A^{-2})(e_1-e_3)+ (A^4-1)(e_2e_1+e_1e_2-e_2e_3-e_3e_2),$$

and thus, 

\begin{equation} \label{penultimateequation}
0 \equiv A^{-2}(A^8-1)(e_1-e_3) +  (A^4-1)(e_2e_1+e_1e_2-e_2e_3-e_3e_2).
\end{equation}

Every Temperley-Lieb element in the equation above intersects the compressing disc $D$ transversely twice. Therefore, we can use Equation (\ref{wid2}) in this situation by carefully taking into account which strings intersect $D$. For example, in the first relation below, the third and fourth strings intersect with $D$. Thus, we get the following equivalences:

\begin{equation}\label{setofsixequations}
\begin{split}
(A^8-1)e_1 & \equiv (A^2-A^6)e_1e_3, \\
(A^8-1)e_3 &  \equiv (A^2-A^6)e_3e_1, \\
(A^8-1)e_1e_2 & \equiv(A^2-A^6)e_3e_1e_2, \\
(A^8-1)e_2e_1 & \equiv(A^2-A^6)e_2e_1e_3, \\
(A^8-1)e_1e_2e_3 & \equiv (A^2-A^6)e_1e_2e_3e_1= (A^2-A^6)e_1e_3, \mbox{and} \\
(A^8-1)e_3e_2e_1 & \equiv (A^2-A^6)e_3e_2e_1e_3= (A^2-A^6)e_3e_1.
\end{split}
\end{equation}

We use the equivalence $(A^8-1)e_1 \equiv (A^2-A^6)e_1e_3\equiv(A^8-1)e_3$ and therefore, the first two terms in Equation (\ref{penultimateequation}) cancel out and we get the following relation:

\begin{equation}\label{finalequation}
0 \equiv (A^4-1)(e_2e_1 -e_2e_3 +e_1e_2-e_3e_2).\end{equation}

We now embed $TL_4$ into $H_3 = F_{0,4} \times I$ as described earlier and under the homomorphism $\rho_*$, we get that $e_1e_2 \mapsto a_1a_3[a_2a_3]$, \  $e_3e_2 \mapsto a_2a_3[a_1a_3]$, $e_2e_1 \mapsto [a_1a_2]$, and $e_2e_3 \mapsto [a_1a_2]$. In particular, $\rho_*(e_2e_1) = \rho_*(e_2e_3)$ as illustrated in Figure \ref{H2H1counterexamplefigure}. Here $[a_ia_j]$ represents a curve that separates the boundary components $a_i$ and $a_j$ from the other two boundary components of $F_{0,4}$.

\begin{figure}[h]
    \centering
$$\vcenter{\hbox{
    \begin{overpic}
    [unit=1mm, scale = 0.55]{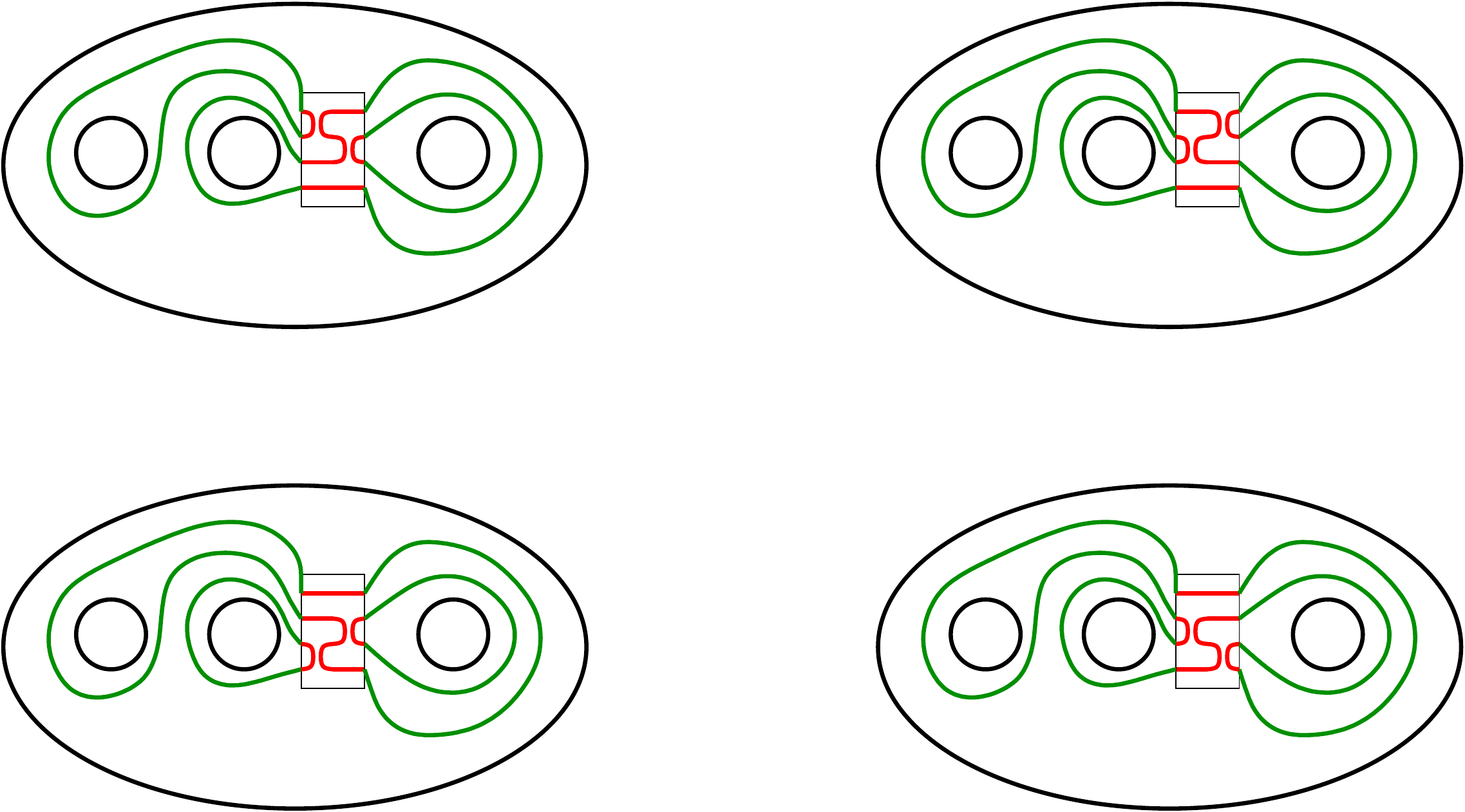}
\put(26,50){$\footnotesize{{e_1e_2}}$}
\put(26,8){$\footnotesize{{e_3e_2}}$}
\put(102,50){$\footnotesize{{e_2e_1}}$}
\put(102,8){$\footnotesize{{e_2e_3}}$}
\put(20,38){$a_1a_3[a_2a_3]$}
\put(20,-4){$a_2a_3[a_1a_3]$}
\put(99,38){$[a_1a_2]$}
\put(99,-4){$[a_1a_2]$}
    \end{overpic}}}$$
    \caption{{\color{white}.}}
    \label{H2H1counterexamplefigure}
\end{figure}

Thus, in $\mathcal{S}_{2,\infty}(H_3)_{\gamma}$, Equation (\ref{finalequation}) results in the following equivalence: 

\begin{equation}\label{newrelation}
0 \equiv (A^4-1)(a_1a_3[a_2a_3] - a_2a_3[a_1a_3]).\end{equation} 

This relation consists of two curve systems that are not  ambient isotopic in $F_{0,4}$. Notice that in Theorem \ref{przytyckitheorem} the ideal $\mathcal{I}$ has the following generators: for every even $k$ and $z_k$ having minimal intersection number $k$ with $D$, the ideal $\mathcal{I}$ has exactly one generator $ (A^{2k+4}-1)z_k + \sum\limits_{i=0}^{k-2} \alpha_i(A)z_i.$ Therefore, the right hand side of Equation (\ref{newrelation}) is not contained in the ideal $\mathcal{I}$ and thus, we have found a new relation in $\mathcal{S}_{2,\infty}(H_3)_{\gamma}$ which serves as a counterexample to Theorem \ref{przytyckitheorem}. This completes the proof of part (1) of Theorem \ref{counterexampletheoremconnsum}. \\

To prove part (2) of Theorem \ref{counterexampletheoremconnsum}, we observe that $H_2 \ \# \ H_1$ can be embedded in the connected sum of any two handlebodies of higher genera and the same curve system in Figure \ref{H2H1counterexamplefigure} embedded in the surface $F_{0, n+m+1}$ leads to a counterexample for all connected sums $H_n \ \# \ H_m$, $n + m \geq  3$.

\end{proof}

\begin{remark}

When we compare the sliding relations $\phi_t$ with $\overline{\phi_t}$ we obtain the equivalence $0 \equiv (A^4 -1)^2(e_1 +e_3-e_1e_2e_3-e_3e_2e_1)$ (see the calculation below). For $H_n \ \# \ H_1$ this relation cannot give a counterexample as it vanishes after embedding the Temperley-Lieb box into $H_n \ \# \ H_1$ since $\rho_*(e_1) = \rho_*(e_3) = \rho_*(e_1e_2e_3) = \rho_*(e_3e_2e_1)$ (see Figure \ref{remarkembeddingtlk}). However, this relation is still nontrivial if we embed the Temperley-Lieb box into $H_n \ \# \ H_2$ as in Figure \ref{z4figure}. We leave this as an exercise to the reader.

\begin{figure}[h]
    \centering
    \includegraphics[scale=0.45]{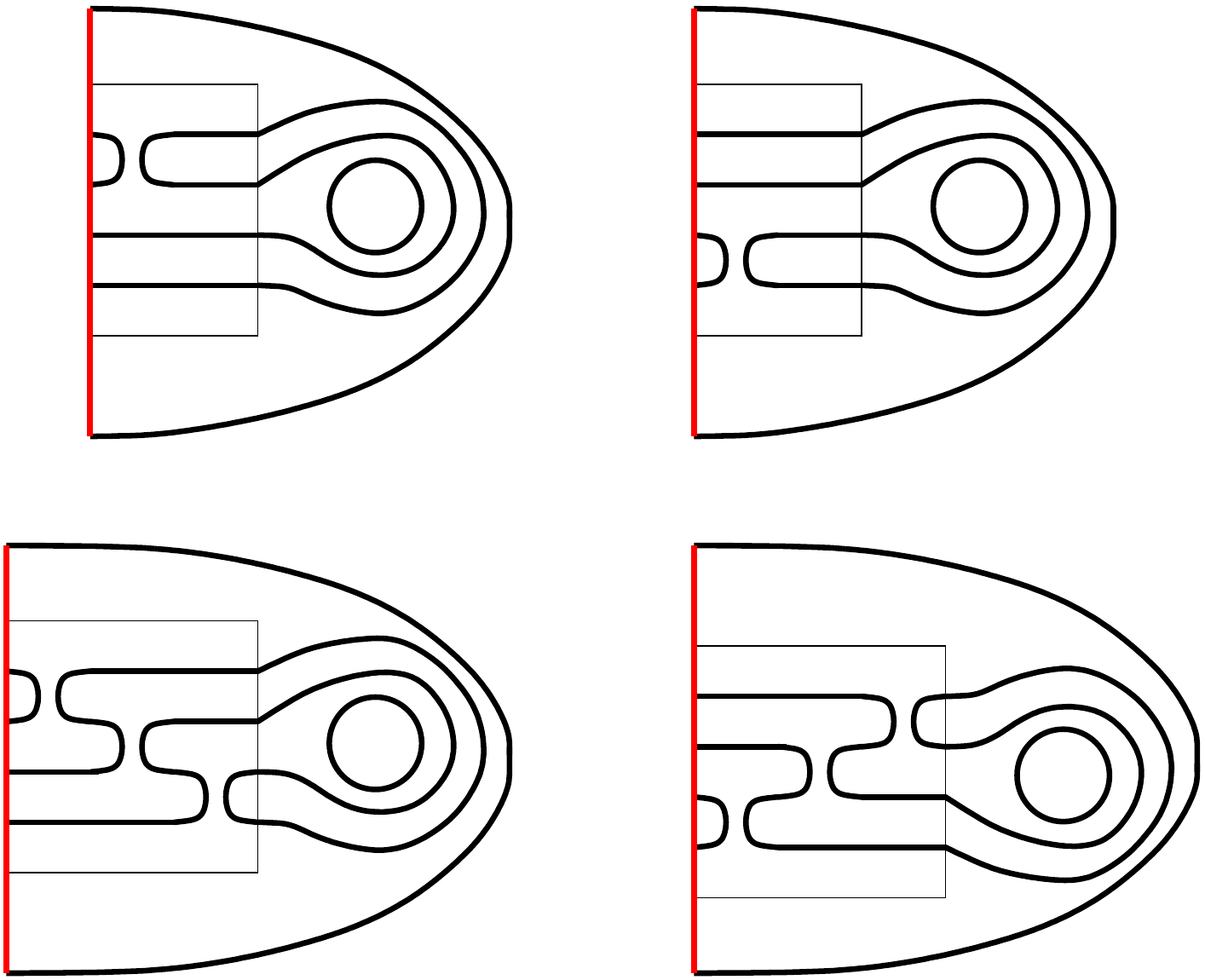}
    \caption{ Embedding Temperley-Lieb elements $e_1$,$e_3$, $e_1e_2e_3$, and $e_3e_2e_1$}
    \label{remarkembeddingtlk}
\end{figure}

Calculation: Consider the sliding relation given by $\phi_t$ in Equation (\ref{connsumeqn2}). Multiplying the sliding relation given by its mirror image $\overline\phi_t$ by $A^{12}$ and adding it to Equation (\ref{connsumeqn2}) we get:
{\color{black}\begin{dmath} 
0 \equiv A^2(A^4-1)^2(e_3-e_3e_2e_1-e_1e_2e_3 - e_3e_1e_2-e_2e_1e_3) + A^{-2}(A^{12}-1)(1-A^4)e_1 + (A^8-1)(1-A^4)(e_1e_3+e_2e_1+e_1e_2).
\end{dmath}}

Now the terms which intersect the compressing disc in exactly two 
points (for example, $e_1$, $e_3$, $e_1e_2$, $e_2e_1$, $e_1e_2e_3$,
and $e_3e_2e_1$) satisfy 
Equation (\ref{setofsixequations}).
Three terms $e_1e_3$, $e_1e_3e_2$, and $e_2e_3e_1$ are disjoint from compressing disc. Thus, after reduction we get the required equivalence:

{\color{black}$$0 \equiv (A^4-1)^2(e_1+e_3-e_1e_2e_3 -e_3e_2e_1).$$}

\end{remark}

\section{Future Directions}

We have shown that Theorem \ref{przytyckitheorem} does not hold in full generality. However, our calculations suggest that the sliding relations that generate the ideal $\mathcal{I}$ are enough in the case of $H_1 \ \# \ H_1$. 
    
\begin{conjecture}\label{h1h1conjecture}

$\mathcal{S}_{2,\infty}(H_1 \ \# \ H_1) = 
\mathcal{S}_{2,\infty}(H_2)/\mathcal{I}. $

\end{conjecture}

In support of the conjecture we have checked that when $k=4$, all the handle sliding relations come from $\phi_t$ and sliding relations from the case $k=2$, and when $k=6$, all the handle sliding relations again come from $\phi_t$ and sliding relations from the smaller cases $k=2$ and $4$. \\

In a future paper we plan to resolve this conjecture and as an application use it to compute the Kauffman bracket skein module of the connected sum of lens spaces over the ring $\mathbb{Z}[A^{\pm 1}]$. In particular, we will compare our result with the result in \cite{rp3rp3} about the connected sum of two copies of the real projective space, $\mathbb{R}P^3 \ \# \ \mathbb{R}P^3$. \\

\section{Acknowledgements}

The second author was partially supported by Simons Collaboration Grant-637794 and the CCAS Enhanced Travel award. The authors would like to thank Charles Frohman  due to whom they decided to provide a proof for Theorem \ref{connsum} and ended up disproving it.\footnote{On March 17, 2020, already during the COVID-19 pandemic, Charles Frohman had asked the second author whether he had published a proof of Theorem \ref{przytyckitheorem}.}

\end{document}